\documentclass[12pt]{article}
\usepackage{fullpage}
\usepackage{amsmath,amsfonts,amssymb,amsthm,graphics,amscd,MnSymbol}
\usepackage{graphicx}
\usepackage{color}

\usepackage{makecell}
\usepackage{multirow}
\usepackage{subcaption}
\usepackage{booktabs}
\usepackage{placeins}

\usepackage{tabularx, booktabs}
\usepackage[
pdfauthor={},
pdftitle={The star edge coloring of cubic Halin graphs with star chromtic index $5$},
pdfstartview=XYZ,
bookmarks=true,
colorlinks=true,
linkcolor=blue,
urlcolor=blue,
citecolor=blue,
bookmarks=false,
linktocpage=true,
hyperindex=true
]{hyperref}

\begin{document}

\newtheorem {Theorem}  {Theorem}[section]
\newtheorem {Lemma}[Theorem]{Lemma}
\newtheorem {Conjecture}[Theorem]{Conjecture}
\newtheorem {Proposition}[Theorem]{Proposition}
\newtheorem {Corollary}[Theorem]{Corollary}
\newtheorem {Problem}[Theorem]{Problem}
\newtheorem {Question}[Theorem]{Question}
\newtheorem {Claim}{Claim}
\newtheorem {Observation}[Theorem]{Observation}
\theoremstyle{definition}
\newtheorem{Definition}[Theorem]{Definition}
\newtheorem{Example}[Theorem]{Example}
\newenvironment {Proof} {\noindent {\bf Proof.}}{\quad $\square$\par\vspace{3mm}}

\def\FR#1#2{\frac{#1}{#2}}
\def\grad{gradual}
\def\NN{{\mathbb N}}
\def\bQ{{\bf Q}}
\def\hS{{\hat S}}
\def\hQ{\hat{\bf Q}}
\def\hu{{\hat u}}
\def\hP{{\hat P}}
\def\hk{{\hat k}}
\def\vareps{\varepsilon}
\def\eps{\varepsilon}
\def\VEC#1#2#3{#1_{#2},\dots,#1_{#3}}
\def\Gb{\overline{G}}
\def\Kb{\overline{K}}
\def\FL#1{\left\lfloor{#1}\right\rfloor}
\def\CL#1{\left\lceil{#1}\right\rceil}
\def\CH#1#2{\binom{#1}{#2}}
\def\esub{\subseteq}
\def\bT{{\mathbf T}}
\def\bB{{\mathbf B}}
\def\bC{{\mathbf C}}
\def\C#1{{\left|#1\right|}}
\def\diam{{\rm diam}}
\def\rad{{\rm rad}}
\long\def\skipit#1{}
\def\gjoin{\diamondplus}
\def\SE#1#2#3{\sum_{#1=#2}^{#3}}
\def\la{\langle}
\def\ra{\rangle}

\date{\today}

\title{The star edge coloring of cubic Halin graphs with star chromatic index $5$}

\author{
Xingxing Hu~\thanks{Department of Mathematics, China Jiliang University, Hangzhou 310018, China. },\,\,
Yunfang Tang~\thanks{
Department of Mathematics, China Jiliang University, Hangzhou 310018, China. tangyunfang8530@163.com. Supported by the National Natural Science Foundation
of China (Grant 11701543).}}

\maketitle

\vspace{-2pc}
\begin{abstract}
The star chromatic index of a graph $G$, denoted by $\chi'_{st}(G) $, is the minimum number of colors needed to properly color the edges of $G$ such that no path or cycle of length four is bi-colored. Casselgren et al. and Hou et al. independently proved that the star chromatic index of a cubic Halin graph, except a special graph, is at most $6$. It remains an open problem  to determine which of such graphs have star chromatic index $5$. In this paper, we show that if $G\ne N_{e_2}$ is a cubic Halin graph whose tree is a caterpillar or a complete tree, then $\chi'_{st}(G)=5$.

\noindent {\bf Keywords:}  star chromatic index; caterpillar; complete tree; Halin graph
\end{abstract}

\baselineskip 16pt

\section{Introduction }
 \hspace*{2em}
 All graphs in this paper are finite and simple. Given a graph $G=V(G)\cup E(G)$, let $c$: $E(G) \rightarrow [k] $ be a proper edge-coloring of $G$, where $k\ge 1$ and $[k] \:=\left\{ 1,2,\cdots k \right\} $. We say that $c$ is a \textit{star $k$-edge coloring} of $G$ if no path or cycle of length four in $G$ is bi-colored under the coloring $c$; and $G$ is \textit{star $k$-edge colorable} if $G$  admits a star $k$-edge coloring. The \textit{star chromatic index} of $G$, denoted by $\chi'_{st}(G)$, is the smallest integer $k$ such that $G$ admits a star $k$-edge coloring.

Star edge-coloring of a graph $G$ was introduced by Liu and Deng~\cite{LD2008}, motivated by the star vertex coloring problem.
Dvo\v{r}\'ak et al~\cite{DMS2013} proposed the following conjecture:
\begin{Conjecture} {\rm\cite{DMS2013}}\label{c1}
If $G$ is a subcubic graph, then $\chi'_{st}(G) \leq 6$.
\end{Conjecture}
Although Conjecture~\ref{c1} remains open, it was confirmed true for some special graphs, such as subcubic outer-planar graphs~\cite{BL2016}, subcubic graphs with maximum average degree at most $5/2$~\cite{lei2018,lei2018b}.
In addition, star edge-coloring of other aspects have also been studied extensively, which was explicitly introduced in a survey~\cite{LS2021}.

A \textit{Halin graph} $G=T \cup C$ is a planar graph that consists of plane embedding of a tree $T$ that has no vertices of degree two, and a cycle $C$ connecting all the leaves of $T$ such that $C$ is the boundary of the exterior face. The tree $T$ is called the \textit{characteristic tree}. If every non-leaf vertex in the tree $T$ of a Halin graph $G$ has degree exactly $3$, then $G$ is called a \textit{cubic Halin graph}.

Recently, a cubic Halin graph was proved true by Casselgren et al~\cite{JCC2021} and Hou et al.~\cite{Hou2019}, independently.

For $l\ge1$, a \textit{complete tree} $T_l$ is a tree of height $l+1$ with a root vertex $v_0$ such that all its leaves are at the same distance $l$ from $v_0$, we call $T_l$ has $l$ levels. A \textit{complete cubic Halin graph} is a cubic Halin graph whose characteristic tree is a complete tree $T_l$, denoted by $G_l$. Clearly, $T_1\cong K_4$. A \textit{caterpillar} is a tree such that the removal of all its leaves results in a path. Let $\mathcal{G}_h$ be the set of all cubic Halin graphs whose characteristic trees are caterpillars with $h+2$ leaves. Specially, $\mathcal{G}_1=\{K_4\}$. For $h\ge 2$, a Halin grpah $G=T\cup C_{h+2}\in \mathcal{G}_h$, denote the spine of $T$ by $P_h=v_1v_2\ldots v_h$, let $u_0,u_1$, and $u_h,u_{h+1}$ be the neighbors of $v_1$ and $v_h$ on $C_{h+2}$, respectively.  For $2\le i\le h-1$ ($h\ge 3$), let $u_i$ be the neighbor of $v_i$ on $C_{h+2}$, then $u_iv_i\in E(G)$, we call such an edge $u_iv_i$ a leaf-edge.  We draw the graph $G$ on the plane such that the spine $P_{h}$ in the middle, and the leaf-edges $u_iv_i$ incident with $v_i$ either on the left side or right side of $P_{h+2}$ for all $i\in \{2,\ldots, h-1\}$. Specially, if all the leaf-edges are on the same side of $P_h$, as we all know, $G$ is called a necklace graph, denoted by $N_{e_h}$.

A strong edge coloring of a graph $G$ is a proper edge coloring so that no edge can be adjacent to two edges with the same color. A strong edge coloring is a special case of star edge coloring, therefore, the strong chromatic index of a graph $G$, denoted by $\chi'_s(G)$, is an upper bound of $\chi'_{st}(G)$. About the results of strong edge coloring of cubic Halin graphs, we will see~\cite{CL2012,LL2012,Shiu2006,ShiuT2009,WY2021}.
For a complete cubic Halin graph $G$, Shiu et al.~\cite{ShiuT2009} determined $\chi'_{st}(G)=6$ except $G\ne N_{e2}$.
For $G\in \mathcal{G}_h$, Shiu et al.~\cite{Shiu2006} determined $\chi'_{st}(N_{eh})$, showed that $6\leq\chi'_{s}(G) \leq 8$ if $h\ge 4$, and conjectured that $\chi'_{st}(G)=6$ if $h$ is odd.
Furthermore, Chang et al.~\cite{CL2012} disproved the conjecture of Shiu et al.and determined the values of star chromatic indices
for some special families of graphs $G\in \mathcal{G}_h$. Up to date, there are very few
results about star chromatic index for the above two kinds of cubic Halin graphs. Hou et al.~\cite{Hou2019} proved the following result.
\begin{Theorem}{\rm\cite{Hou2019}}\label{Th0}
If $h$ is odd, then $4 \leq \chi'_{st}(N_{eh}) \leq 5$.
\end{Theorem}
Hence, based on the above known results, it is interesting to study star edge coloring of the above two kinds of cubic Halin graphs. In this paper, we determine the exact values of star chromatic indices for them,  which improves Theorem~\ref{Th0}, and extends the known results. Specifically, in Section $1$, we proved that $\chi _{st}^{'}(G) \geq5$ for any cubic graph $G$ (see Theorem \ref{key3}); in Section $2$, we showed that if $G\ne N_{e_2}$ is a complete cubic Halin graph, then $\chi _{st}^{'}(G)=5$ (see Theorem~\ref{Th1}); in Section $3$, we showed that if $G\in \mathcal{G}_h$ and $G\ne N_{e_2}$, then $\chi _{st}^{'}(G)=5$ (see Theorem~\ref{Th2}), and hence $\chi _{st}^{'}(N_{e_h})=5$ for $h\ne 2$ (see Corollary~\ref{c1}).

\section{Prelimilaries}
\begin{Theorem} \label{key3}
If $G$ is a cubic Halin graph, then $\chi'_{st}(G) \geq5$.
\end{Theorem}

\begin{figure}[h]
 \captionsetup{labelsep=period}
 \centering
 \begin{subfigure}[t]{0.3\linewidth}
  \centering
  \includegraphics[width=0.6\linewidth]{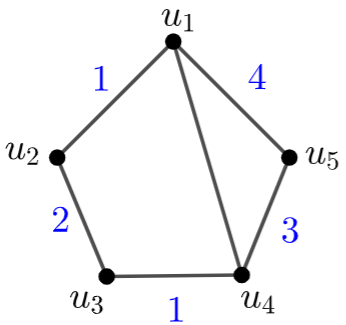}
  \caption{$u_1u_4\in E(T)$}
  \label{fig:ab1-a}
 \end{subfigure}
 \hspace{0.2\linewidth}
 \begin{subfigure}[t]{0.3\linewidth}
  \centering
  \includegraphics[width=0.6\linewidth]{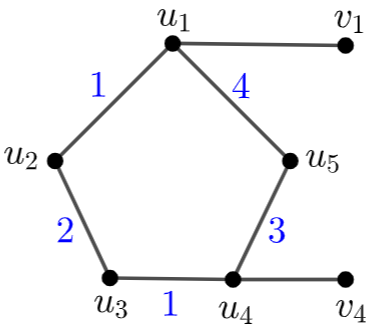}
  \caption{$u_1u_4\notin E(T)$}
  \label{fig:ab1-b}
 \end{subfigure}
 \caption{a-b}
 \label{fig:a-b}
\end{figure}

\begin{Proof}
Note that $\chi'_{st}(K_4)=5$, then we may assume $G \not\cong K_4$. We prove the result by contradiction.
Assume $\varphi$ is a star $4$-edge coloring of $G$ using the color set $\{1, 2, 3, 4\}$.
It is easy to see $G$ contains a $5$-cycle $C_5$. Let $C_5=u_1u_2u_3u_4u_5u_1$. Note that $C_5$ is star $4$-edge colorable. Therefore, there exist two non-adjacent edges in $C_5$ with the same color. Without loss of generality, we may assume $\varphi(u_1u_2) = \varphi(u_3u_4)=1$, $\varphi(u_2u_3)=2$, $\varphi(u_4u_5)\}=3$, and $\varphi(u_5u_1)=4$. Note that $d_G(u_1)=d_G(u_4)=3$.
Then $u_1u_4\in E(G)$ or $u_1,u_4$ has another neighbor $v_1,v_4$  not on  $C_5$, respectively (see Figure~\ref{fig:a-b}).

If $u_1u_4\in E(G)$, then $\varphi(u_1u_4)\in [4]\setminus \{1,3,4\}=\{2\}$, since $u_1u_4$ is adjacent to the four edges of $C_5$ except $u_2u_3$. But $\varphi(u_1u_4)=\varphi(u_2u_3)$, which results in a bi-chromatic $4$-cycle $u_1u_2u_3u_4u_1$ with $1-2-1-2$, a contradiction.

Suppose that $u_1v_1,u_4v_4\in E(G)$. Then $\varphi(u_1v_1)\in \{2,3\}$, $\varphi(u_4v_4)\in \{2,4\}$.
We claim that $2\notin\{\varphi(u_1v_1),\varphi(u_4v_4)\}$. Otherwise, there is a bi-chromatic $4$-path $v_1u_1u_2u_3u_4$
or $u_1u_2u_3u_4v_4$ with $2-1-2-1$ or $1-2-1-2$, correspondingly. Then $\varphi(u_1v_1)=3$, and $\varphi(u_4v_4)=4$, which implies in a
bi-chromatic $4$-path $v_1u_1u_5u_4v_4$ with $3-4-3-4$, a contradiction.

From the above, the result follows.
\hfill\end{Proof}

\section{Characteristic trees are complete trees}
\hspace{2em}
Let $T_0^i$ be a full binary tree of 3 levels, where the root vertex is $v_0^i$, the first-level vertices are $x^i, y^i$, the second-level vertices are sequentially $x_1^i, x_2^i, y_1^i, y_2^i$, and the third-level vertices are sequentially $x_{11}^i, x_{12}^i, x_{21}^i, x_{22}^i, y_{11}^i, y_{12}^i, y_{21}^i, y_{22}^i$, which are leaves. Define $F^i$ as the graph obtained from $T_0^i$ by sequentially connecting the leaves of $T^i$, such that there is a path $P_8 = x_{11}^i x_{12}^i x_{21}^i x_{22}^i y_{11}^i y_{12}^i y_{21}^i y_{22}^i$ in $F^i$ (will be applied in Lemma~\ref{key4}). For $i \geq 0$, let $G^i=T^i\cup C^i$ be a cubic Halin graph. For a vertex $v^i\in C^i$, let $u^i$, $s^i$, and $t^i$ be the neighbors of $v^i$, where $u^i\in T^i$, $s^i,t^i\in C^i$. First we delete the two edges  $v^is^i$ and $v^it^i$ in $G^i$; then attach $F^i$ at $v^i$ to $u^i$ such that $v_0^i = v^i$; lastly, add two new edges $x_{11}^i s^i$ and $y_{22}^i t^i$ in $G^i$. Let $G^{i+1}$ be a graph obtained from $G^i$ by the above three steps.
It can be shown that if $G^0 = G_l$, then $G_{l+3}=G^k$, where $k=3 \times 2^{l-1}$.

\begin{Lemma} \label{key4}
For $i \geq 0$, let $G^i$ be a cubic Halin graph. If $\chi'_{st}(G^i) \leq 5$, then $\chi'_{st}(G^{i+1}) \leq 5$.
\end{Lemma}

\begin{figure}[htbp]
    \captionsetup{labelsep=period}
    \centering  % 确保内容居中
    \includegraphics[width=0.56\textwidth]{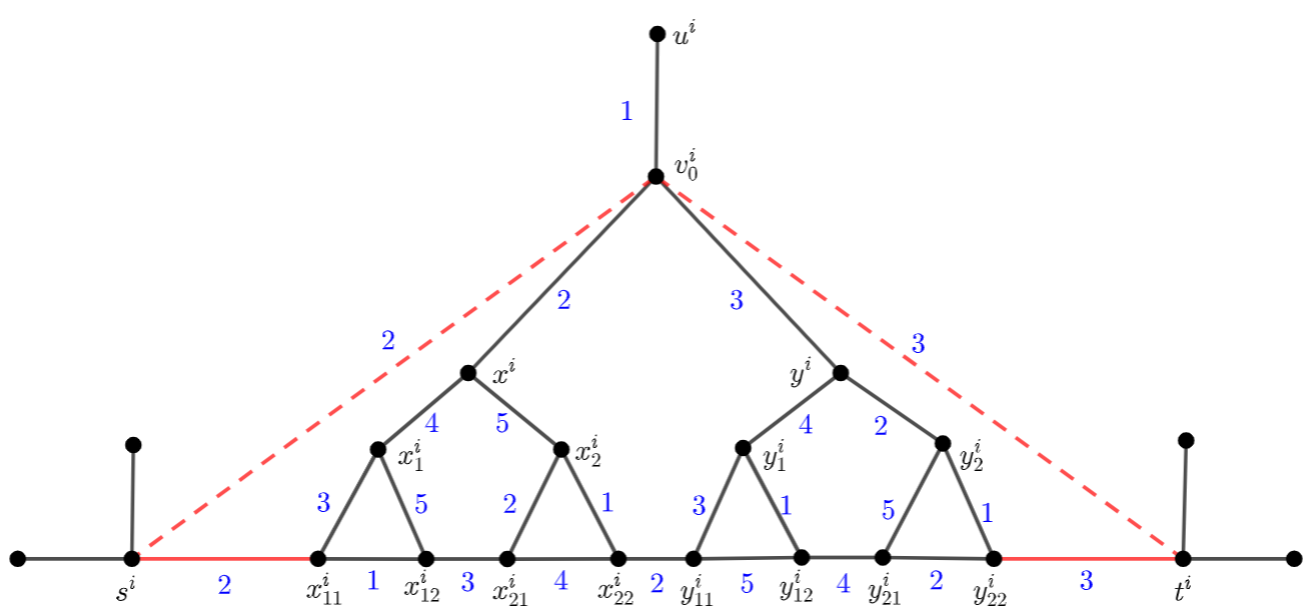}  % 调整宽度以适应页面
    \caption{Structure of $G^{i+1}$.}  % 图片标题
    \label{fig:gi_plus_1}  % 标签用于引用
\end{figure}

\begin{Proof}
Let $\varphi$ be the star $5$-edge coloring of $G^i$. Without loss of generality, assume that $\varphi(v_0^i u_1^i) = 1$,  $\varphi(v_0^i s_1^i) = 2$,  $\varphi(v_0^i t_1^i) = 3$. Now we extend the coloring $\phi$ of $G^i$ to $G^{i+1}$ using the colors shown in Figure~\ref{fig:gi_plus_1}. This gives a star $5$-edge coloring of $G^{i+1}$. Thus the result follows.
\hfill\end{Proof}

\begin{Theorem}\label{Th1}
If $G$ is a complete cubic Halin graph, then $\chi'_{st}(G)=5$.
\end{Theorem}
\begin{figure}[htbp]
    \captionsetup{labelsep=period}
    \centering
    \begin{subfigure}[t]{0.3\linewidth}
        \centering
        \includegraphics[width=0.75\linewidth]{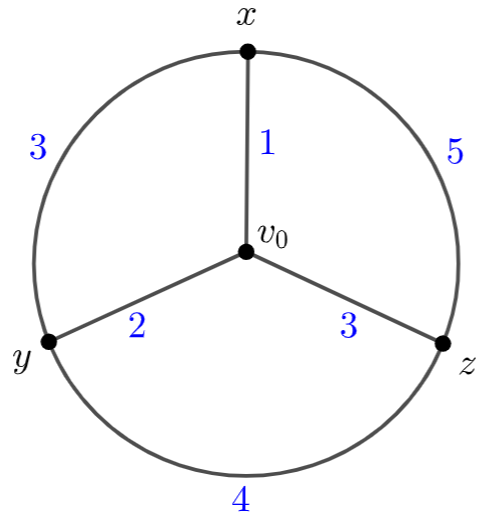}
        \caption{$G_1=W_3$}
        \label{fig3:ab-a}
    \end{subfigure}
    \hfill % Add horizontal space between subfigures
    \begin{subfigure}[t]{0.3\linewidth}
        \centering
        \includegraphics[width=0.8\linewidth]{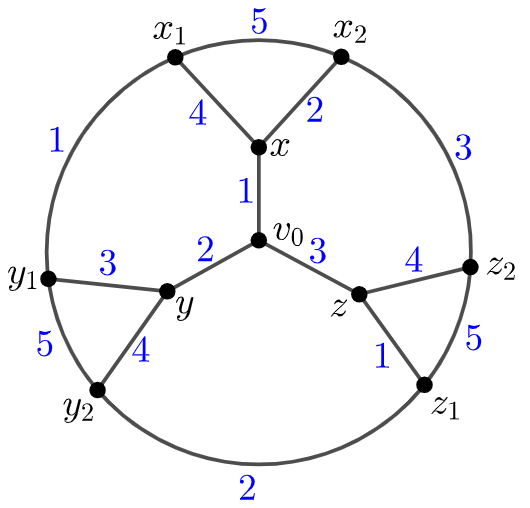}
        \caption{$G_2$}
        \label{fig3:ab-b}
    \end{subfigure}
    \hfill % Add horizontal space between subfigures
    \begin{subfigure}[t]{0.3\linewidth}
        \centering
        \includegraphics[width=0.85\linewidth]{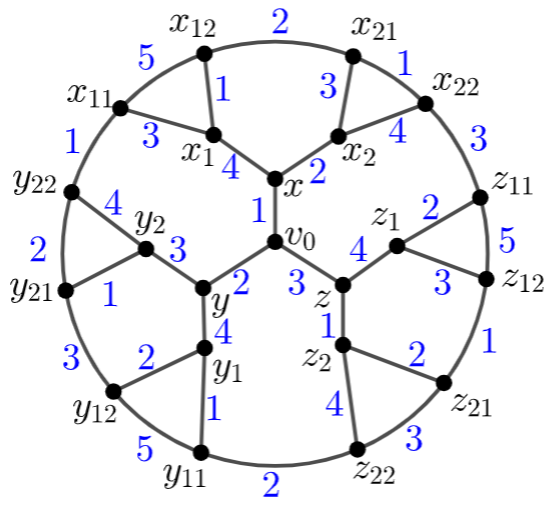}
        \caption{$G_3$}
        \label{fig3:ab-c}
    \end{subfigure}
    \caption{$G_l,l=1,2,3$}
    \label{fig3:abc}
\end{figure}

\begin{Proof}
Let $G=G_l$. If $l=1$, then $G=K_4$. From Figure~\ref{fig3:abc}(\subref{fig3:ab-a}), it is clear that $\chi'_{st}(G)=5$.
Suppose that $l\ge 2$. By Theorem~\ref{key3}, we only need to prove that $\chi'_{st}(G) \leq 5$.
We will prove it by induction on the level $l$ of $G$.
If $l=2, 3$, then  $\chi'_{st}(G) \leq 5$, see Figure~\ref{fig3:abc}(\subref{fig3:ab-b}) and ~\ref{fig3:abc}(\subref{fig3:ab-c}). Suppose that for any graph $G$ with level $l \leq k$, the result holds. Repeatedly using Lemma~\ref{key4}, we can directly prove that if $l =k+3$, then $\chi'_{st}(G) \leq 5$. Thus the result follows.
\hfill\end{Proof}

\section{Characteristic trees are caterpillars}
\begin{Theorem}\label{Th2}
 Let $G=T\cup C\in \mathcal{G}_h$ and $G\ne Ne_2 $, then $\chi'_{st}(G)=5$.
\end{Theorem}

\begin{Proof}
Note that $\mathcal{G}_h=\left\{ Ne_h \right\}$ for $h=1,2,3$. It is easy to check that $\chi'_{st}( Ne_1 ) =\chi'_{st}(Ne_3) =5$ and $\chi'_{st}(Ne_2) =6$. Suppose that $h\ge 4$. By Theorem~\ref{key3}, we only need to prove $\chi'_{st}(G) \le 5$.
We prove it by induction on $h$. For $h=4,5$, $\chi'_{st}(G) \le 5$, all such colorings are supplied in Figure~\ref{fig:a-e}.
\begin{figure}[h]
 \captionsetup{labelsep=period}
 \centering
 \begin{subfigure}[t]{0.3\linewidth}
  \centering
  \includegraphics[width=0.8\linewidth]{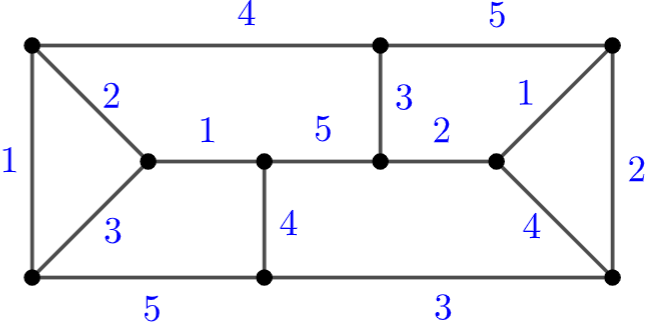}
  \caption{}
  \label{fig:ab2-a}
 \end{subfigure}
 \hspace{0.08\linewidth}
 \begin{subfigure}[t]{0.3\linewidth}
  \centering
  \includegraphics[width=0.8\linewidth]{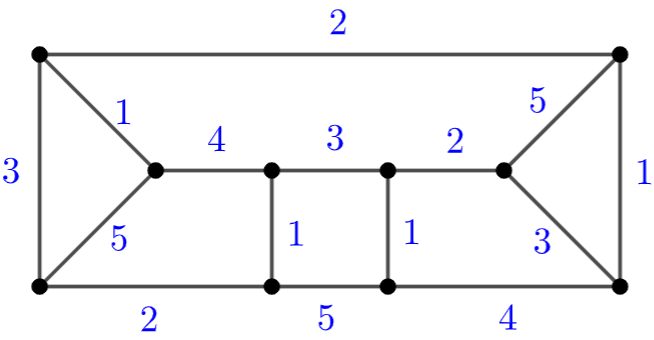}
  \caption{}
  \label{fig:ab2-b}
 \end{subfigure}
 \begin{subfigure}[t]{0.3\linewidth}
  \centering
  \includegraphics[width=0.8\linewidth]{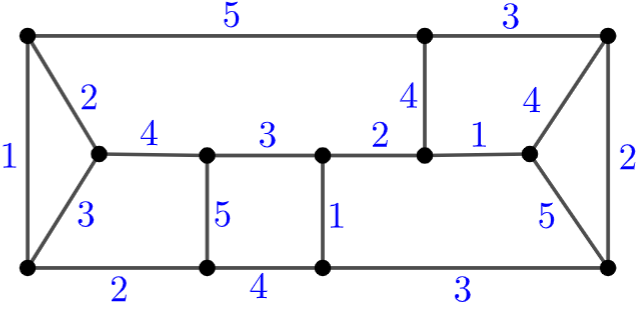}
  \caption{}
  \label{fig:ab2-c}
 \end{subfigure}
 \begin{subfigure}[t]{0.3\linewidth}
  \centering
  \includegraphics[width=0.8\linewidth]{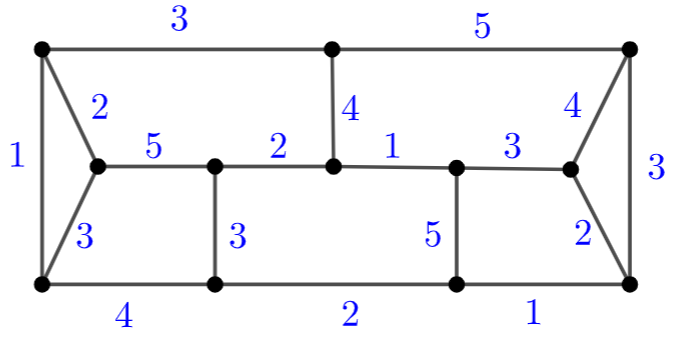}
  \caption{}
  \label{fig:ab2-d}
 \end{subfigure}
 \begin{subfigure}[t]{0.3\linewidth}
  \centering
  \includegraphics[width=0.8\linewidth]{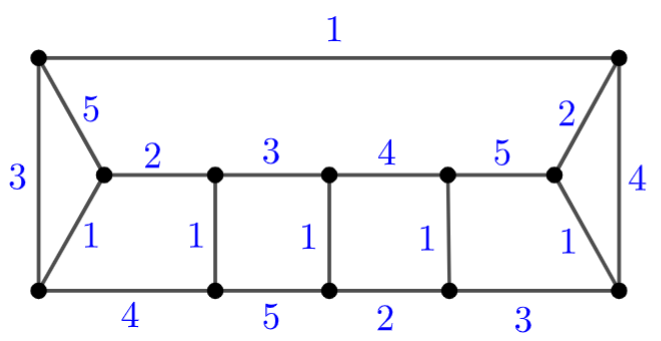}
  \caption{}
  \label{fig:ab2-e}
 \end{subfigure}
 \caption{a-e}
 \label{fig:a-e}
\end{figure}

In the following, we may assume $h\geq6$.
Let $P=v_1v_2\cdots v_{h+2}$ be the longest path in $T$. Rename the vertices so that $u_1=v_1$, $v=v_2$, $w=v_3$,$x=v_4$ and $y=v_5$. Except $u_1$ and $v$, let the other neighbor of $u$ is a leaf $u_2$ on $C$. Except $u,u_1$ ($u,u_2$), let $x_1$ ($y_1$) be the other neighbor of $u_1$ ($u_2$). Except the edge $x_{i-1}x_i$ ($y_{i-1}y_i$), let the other two  edges incident to $x_i$ ($y_i$) in $G'$ be $x_ix'_i$ ($y_iy'_i$) and $x_ix_{i+1}$ ($y_{j-1}y_j$) for $i\in\{1,2,3\}$ (let $x_0=u_1$, $y_0=u_2$). For convenience, let $x_i,y_i\in V(C)$, and $x'_i,y'_i\in V(P)$.

 Note that $x'_1=v$ or $y'_1=v$. Without loss of generality, we may assume $y'_1=v$, i.e. $vy_1\in E(T)$.
In our later inductive steps, let $G'=G-u-u_1-u_2-y_1+vx_1+vy_2$. Note that $G'\in\mathcal{G}_{h-2}$. Then by the induction hypothesis, there exists a star $5$-edge coloring $f'$ for $E(G')$.  For each edge $e$ in $E(G^{'})\backslash\{vx_1,vy_2\}$, let $f(e)=f'(e)$. Without loss of generality, assume that $f(vw)=1$, $f(vx_1)=2$, and $f(vy_2)=3$. Denote by $C(v)=\{f(e):e\in E(v)\}$.
For convenience, let $\{t'_i\}=\{4,5\}\setminus \{t_i\}$ if $t_i\in \{4,5\}$ for $i\in\{0,1,2\}$;
$\{\alpha'_0\}=\{1,2\}\setminus \{\alpha_0\}$ if $\alpha_0\in \{1,2\}$.

\begin{figure}[htbp]
 \captionsetup{labelsep=period}
 \centering
 \begin{subfigure}[t]{0.35\linewidth}
  \centering
  \includegraphics[width=\linewidth]{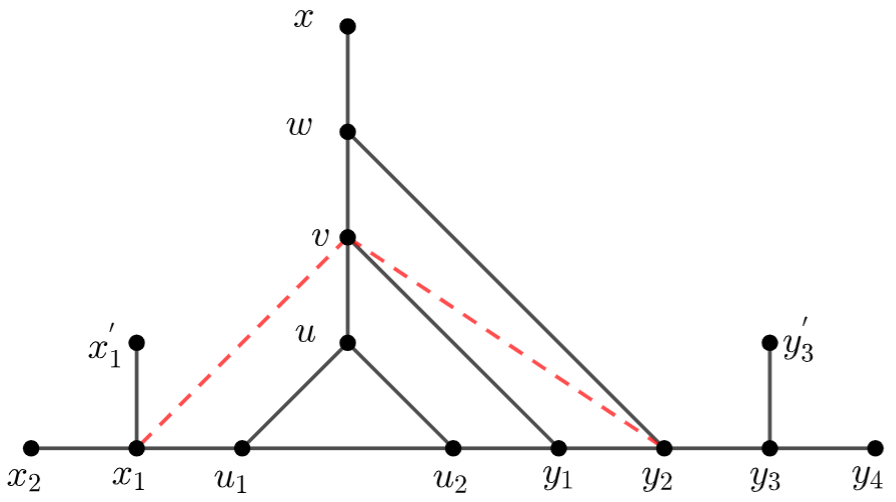}
  \caption{$wy_2\in E\left( T \right)$}
  \label{fig5:ab3-a}
 \end{subfigure}
 \hspace{0.1\linewidth}
 \begin{subfigure}[t]{0.35\linewidth}
  \centering
  \includegraphics[width=1.2\linewidth]{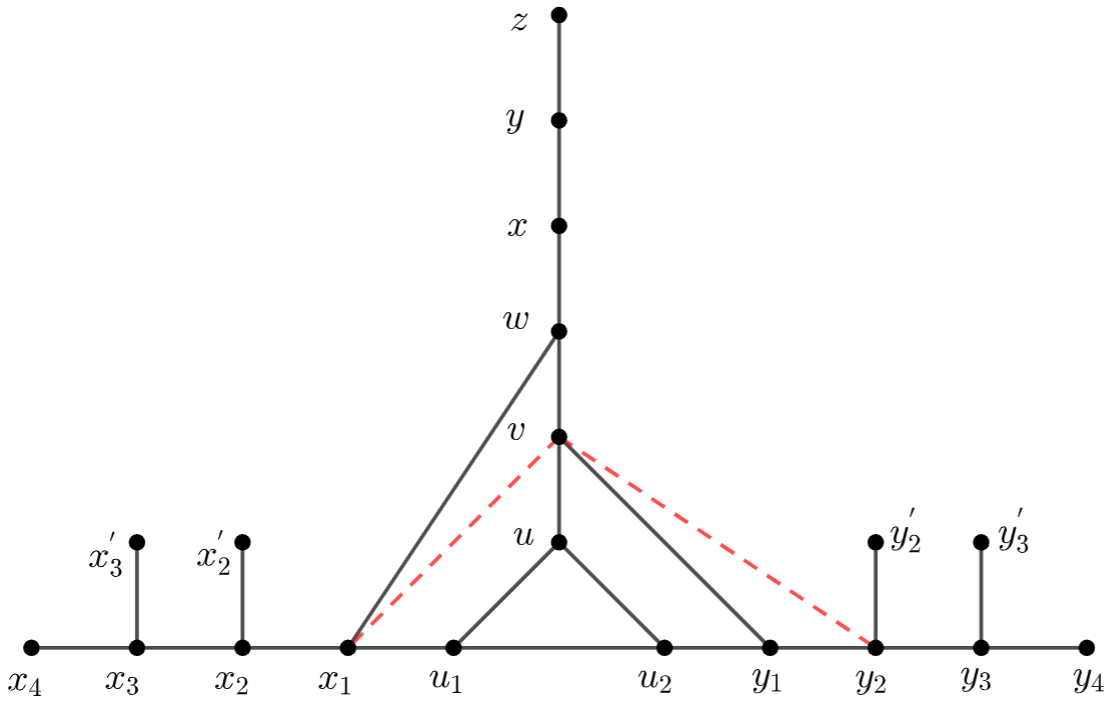}
  \caption{$wx_1\in E\left( T \right)$}
  \label{fig5:ab3-b}
 \end{subfigure}
 \caption{a-b}
 \label{fig5:ab}
\end{figure}

In the following, we will consider two cases to extend $f'$ of $G'$ to the remaining edges of $G$ to get a star  $5$-edge coloring $f$ of $G$.

\textbf{Case1} $w=y'_2$, i.e. $wy_2\in E(T)$ (see Figure~\ref{fig5:ab}(\subref{fig5:ab3-a})).

For convenience, let $f(xw)=t_0$, $f(y_2w)=t_1$, $f(y_2y_3)=t_2$, $f(x_1x'_1)=\lambda_1$, $f(x_1x_2)=\lambda'_1$, $f(y_3y'_3)=\mu_1$ and $f(y_3y_4)=\mu'_1$. We first let $f(u_1x_1)=2$.

\textbf{Subcase 1.1}  $t_0=2$.

Then $t_1\in\{4,5\}$, $t_2\in\{1,2,t'_1\}$.
Let $f(y_1y_2)=f(uu_1)=3$, $f(u_1u_2)=1$, $f(vy_1)=2$, $f(uu_2)\in\{t_1,t'_1\}\setminus \{f(uv)\}$, and
\begin{eqnarray*}
f(uv)=f(u_2y_1)&=&\begin{cases}
t_1& \mbox{if $t_2=t'_1$,} \\
t'_1&\mbox{if $t_2\in\{1,2\}$.}
\end{cases}
\end{eqnarray*}

\textbf{Subcase 1.2} $t_0=3$.

Then $t_1\in\{4,5\}$, $t_2\in\{2,4,5\}$.

Suppose that $t_1\in\{4,5\}$.
If $t_2=2$, or $t_2=t'_1$ and $\{\mu_1,\mu'_1\}\ne \{1,3\}$, then let $f(uv)=f(u_2y_1)=t'_1$, $f(uu_2)=t_1$, $f(vy_1)=2$, $f(u_1u_2)\in\{1,3\}\setminus\{f(uu_1)\}$, and
\begin{eqnarray*}
f(uu_1)=f(y_1y_2)&=&\begin{cases}
1& \mbox{if $t_2=t'_1$, and $1\notin \{\mu_1,\mu'_1\}$} \\
3&\mbox{otherwise.}
\end{cases}
\end{eqnarray*}
Otherwise, $t_2=t'_1$ and $\{\mu_1,\mu'_1\}=\{1,3\}$.  Note that $C(x)\setminus\{t_0\}=\{2,t'_1\}$. Then $xx_1 \in E(T)$, $f(xy) = 2$, $f(xx_1)=\lambda_1=t'_1$, and $\lambda'_1 \in \{1, 3, t_1\}$. First let $f(uu_1) = 3$. If $\lambda'_1 \in \{1, 3\}$, then let $f(u_2y_1) = 1$, $f(vy_1) = 2$, $f(y_1y_2) = 3$, $f(uv) = f(u_1u_2) = t_1$, and $f(uu_2) = t'_1$. If $\lambda'_1= t_1$, then recolor the edge $wv$ by setting $f(wv) = t'_1$, and let $f(u_1u_2)= f(vy_1) = 1$, $f(uv) = f(y_1y_2)=2$, $f(u_2y_1)=t_1$, and $f(uu_2)=t'_1$.

\textbf{Subcase 1.3} $t_0\in\{4,5\}$.

Then $t_1\in \{2,4,5\}$ and $t_2\in\{1,2,4,5\}$.

\textbf{Subcase 1.3.1} $t_2=1$.

Then $t_1=t'_0$.

If $xx_1\in E(T)$. Then $\{\lambda_1,\lambda'_1\}\ne\{4,5\}$.
Let $f(uv)=f(y_1y_2)=3$, $f(u_2y_1)=f(uu_1)=1$, $f(vy_1)=2$,  $f(u_1u_2)\in\{4,5\}\setminus \{\lambda_1,\lambda'_1\}$, and $f(uu_2)\in\{4,5\}\setminus\{f(u_1u_2)\}$.

Now assume that $xy_3\in E(T)$. First let
$f(uv)=f(u_2y_1)=t_0$, $f(y_1y_2)=3$, $f(uu_2)=t_1$, and $f(vy_1)=2$.
If  $s_0\ne 1$, then let $f(uu_1)=3$, and $f(u_1u_2)=1$.
Suppose that $s_0=1$, then $\{\mu_1,\mu'_1\}=\{2,t_0\}$. In this case, we recolor the edge $wv$ in $G$ by letting $f(wv)=3$, and let
$f(uu_1)=1$, and $f(u_1u_2)=3$.

\textbf{Subcase 1.3.2} $t_2\ne 1$.

Then $t_2=t_0$, and $t_1\in\{2,t'_0\}$; or $\{t_1,t_2\}=\{2,t'_0\}$.
First let $f(uv)=t'_0$, $f(y_1y_2)=f(uu_1)=3$, $f(u_1u_2)=1$ and $f(vy_1)=2$.
If $t_2=t_0$, and $t_1\in\{2,t'_0\}$, then let $f(u_2y_1)=t'_0$, and $f(uu_2)=t_0$.
If $\{t_1,t_2\}=\{2,t'_0\}$, then let $f(u_2y_1)=t_0$, $f(uu_2)=2$.\\

\textbf{Case 2} $w=x'_1$, i.e. $wx_1\in E(T) $. (see Figure~\ref{fig5:ab}(\subref{fig5:ab3-b})).

For convenience, let $f(xy)=s_0$, $f(yz)=s_1$, $f(xw)=t_0$, $f(x_1w)=t_1$, $f(x_1x_2)=t_2$,  $f(x_ix'_i)=\lambda_{i-1}$, $f(x_ix_{i+1})=\lambda'_{i-1}$, $f(y_iy'_i)=\mu_{i-1}$ and $f(y_iy_{i+1})=\mu'_{i-1}$, where $i\in \{2,3\}$.

\textbf{Subcase 2.1}  $t_0=3$.

Then $t_1\in \{4,5\}$, and $t_2\in \{1,3,4,5\}$. Let $f(u_2y_1)=1$, $f(u_1x_1)=f(vy_1)=2$, $f(y_1y_2)=f(uu_1)=3$,
$f(uu_2)=\{4,5\}\setminus \{f(uv)\}$, and
\begin{eqnarray*}
f(uv)=f(u_1u_2)&=&\begin{cases}
t'_1& \mbox{if $t_2=1$,} \\
t_1&\mbox{if $t_2\ne 1$.}
\end{cases}
\end{eqnarray*}

\textbf{Subcase 2.2}  $t_0\in\{4,5\}$.

Then $t_1=t'_0$, $t_2\in\{t_0,1,3\}$; or $t_1=3$, $t_2\in\{t_0,t'_0\}$.

\textbf{Subcase 2.2.1} $t_1=t'_0$ and $t_2=1$.

Let $\alpha_0\in \{1,2\}$. If $\alpha_0\notin C(x)$, then let $f(u_1x_1)=2$, and $f(u_1u_2)=f(vu)=t_0$, $f(uu_2)=t'_0$,
$f(u_2y_1)=f(vw)=\alpha_0$, $f(vy_1)=\alpha'_0$, and $f(uu_1)=f(y_1y_2)=3$.

Otherwise, $C(x)\setminus\{t_0\}=\{1,2\}$, which implies $xy_2\in E(T)$.
Then  $\mu_1\in \{1,2\}$, and $\mu'_1\in \{1,2,4,5\}$. If $\mu'_1\in\{4,5\}$, then let $f(y_1y_2)=f(vw)=f(uu_1)=3$, $f(uu_2)\in \{4,5\}\setminus\{\mu'_1\}$, $f(u_1u_2)=\mu'_1$,
$f(vu)=f(u_2y_1)=1$, and $f(u_1x_1)=f(vy_1)=2$.
Suppose that $\{\mu_1,\mu'_1\}=\{1,2\}$, then $\{\lambda_1,\lambda'_1\}=\{3,t_0\}$.
We can recolor the edge $wx_1$ by setting $f(wx_1)=2$, and
let $f(uu_1)=f(u_2y_1)=2$, $f(u_1u_2)=1$, $f(uv)=f(y_1y_2)=3$, $f(uu_2)=t_0$, and $f(u_1x_1)=f(vy_1)=t'_0$.

\textbf{Subcase 2.2.2} $t_1=t'_0$, $t_2\in\{t_0,3\}$; or $t_1=3$, $t_2\in\{t_0,t'_0\}$.

Note that $t_2\ne 1$, then first let $f(u_1u_2)=f(vu)=t_1$, $f(u_2y_1)=1$, $f(u_1x_1)=f(vy_1)=2$, and $f(y_1y_2)=3$.
If $t_1=t'_0$, $t_2\in\{t_0,3\}$, then let $f(uu_1)=3$, and $f(uu_2)=t_0$. Otherwise, then let $f(uu_1)=t'_2$, and $f(uu_2)=t_2$.

\textbf{Subcase 2.3} $t_0=2$.

Then $t_1\in\{4,5\}$, $t_2\in \{3,4,5\}$, and $C(x)\setminus\{t_0\}=\{3,t'_1\}$.

\textbf{Subcase 2.3.1} $t_1\in\{4,5\}$, and $t_2=3$.

Note that $C(x)\setminus\{2\}=\{3,t'_1\}$ and $f(x_1x_2)=f(vy_2)=3$. Then $f(xx_2)=\lambda_1=t'_1$ and $\lambda'_1\ne 2$ if $xx_2\in E(T)$;
$f(xy_2)=\mu_1=t'_1$ and $\mu'_1\ne 2$ if $xy_2\in E(T)$. Thus we can recolor the edge $vw$ such that $f(wv)=t'_1$, and let $f(u_1u_2)=t'_1$, $f(uv)=t_1$, $f(uu_2)=f(vy_1)=1$, $f(x_1u_1)=f(u_2y_1)=2$, and $f(uu_1)=f(y_1y_2)=3$.

\textbf{Subcase 2.3.2} $t_1\in\{4,5\}$, and $t_2=t'_1$.

Then $\{\mu_1,\mu'_1\}\subseteq\{1,2,4,5\}$, and $t'_1\notin \{\lambda_1,\lambda'_1\}$.

\textbf{Subcase 2.3.2.1} $\{\mu_1,\mu'_1\}\neq\{4,5\}$.

Let $f(u_1u_2)=1$, $f(x_1u_1)=f(uu_2)=f(vy_1)=2$, $f(uu_1)=f(y_1y_2)=3$, $f(uv)\in \{4,5\}\setminus\{\mu_1,\mu'_1\}$, $f(u_2y_1)\in\{4,5\}\setminus\{f(uv)\}$.

\textbf{Subcase 2.3.2.2} $\{\mu_1,\mu'_1\}=\{4,5\}$.

Then $\{\mu_1,\mu'_1\}=\{t_1,t'_1\}$. For convenience, there are two remarks, one is that if there are some blue edges and blue numbers in $G'$, see Figure~\ref{fig:ab1}--\ref{fig:abc4}, then it implies that we have recolored those edges in $G$, explicitly, see Figure~\ref{fig:abc2}(\subref{fig:ab2-b}),~\ref{fig:abc2}(\subref{fig:ab2-c}), and~\ref{fig:abc4}; the other is that blue bold faced numbers in Figure~\ref{fig:ab1}--\ref{fig:abc4} means the colors of the remaining eight edges for $E(G)\setminus E(G')$.

First we consider the values of $\mu_2$ and $\mu'_2$.

If $3 \notin \{\mu_2, \mu'_2\}$, then the coloring is illustrated in Figure~\ref{fig:ab1}(\subref{fig:ab1-a}).
Suppose that $1 \notin \{\mu_2, \mu'_2\}$. If $\mu'_1 = t_1$, then let the coloring be shown in Figure~\ref{fig:abc2}(\subref{fig:ab2-a}).
Otherwise, $\mu'_1 = t'_1$, then  $xx_2\in E(T)$, which implies that $s_0=t'_1$ and $\lambda_1=3$. If  $\lambda'_1 = 1$, then see the coloring in Figure~\ref{fig:abc2}(\subref{fig:ab2-b}). If $\lambda'_1 \neq 1$, then see the coloring in Figure~\ref{fig:abc2}(\subref{fig:ab2-c}).

Thus we only need to consider the case when $\{\mu_2,\mu'_2\} = \{1,3\}$. If $3 \notin \{\lambda_1, \lambda'_1\}$, the configuration is shown in Figure~\ref{fig:ab1}(\subref{fig:ab1-b}). Suppose that $1 \notin \{\lambda_1, \lambda'_1\}$. If $\mu'_1 = t_1$,  then the remaining coloring was shown in Figure~\ref{fig:ab3}(\subref{fig:ab3-a}). Otherwise, $\mu'_1 = t'_1$, then the coloring is given in Figure~\ref{fig:ab3}(\subref{fig:ab3-b}).

In the following, we only need to consider the case when $\{\lambda_1,\lambda'_1\} = \{\mu_2,\mu'_2\} = \{1,3\}$. Then $2,t_1 \notin \{\lambda_1,\lambda'_1,\mu_2,\mu'_2\}$. Since $s_0 \in \{3, t'_1\}$, we split into two cases.

First, assume $s_0 = 3$. Then $xy_2 \in E(T)$, so $f(xy_2) = \mu_1 = t'_1$ and $\mu'_1 = t_1$. Since $t_1 \notin \{\lambda_1, \lambda'_1\}$ and $2 \notin \{\mu'_1, \mu_2, \mu'_2\}$, we may swap the colors of $xy_2$ and $wx$ in $G'$,  i.e. $f(xy_2)=2$ and $f(wx)=t_0=t'_1$. Then $G'$ also has a star $5$-edge coloring, and it has already been done in Subcase~2.2.2.

Now suppose $s_0 = t'_1$. Then $xx_2 \in E(T)$, thus $\lambda_1 = 3$ and $\lambda'_1 = 1$. We now consider the values of $\{\lambda_2, \lambda'_2\}$. First consider that $t'_1 \notin \{\lambda_2, \lambda'_2\}$. If $\mu'_1 = t_1$, also see the coloring in Figure~\ref{fig:ab3}(\subref{fig:ab3-a}). Otherwise, $\mu'_1 = t'_1$, then the configuration is  also given in Figure~\ref{fig:ab3}(\subref{fig:ab3-b}). If $\{\lambda_2, \lambda'_2\} = \{2, t'_1\}$, then we may recolor the two edges $wx_1,x_1x_2$ such that $f(wx_1)=3$ and $f(x_1x_2)=t_1$, and extend the coloring in $G$, which was shown in Figure~\ref{fig:abc4}(\subref{fig:ab4-a}). If  $2 \notin \{\lambda_2, \lambda'_2\}$, then we have $s_1\in \{1,2\}$, whether $yx_3\in E(T)$ or $yy_2\in E(T)$, explicitly, see the coloring in Figure~\ref{fig:abc4}(\subref{fig:ab4-b}) and Figure~\ref{fig:abc4}(\subref{fig:ab4-c}).

From the above, it is easy to see that $f$ is a star $5$-edge coloring for $G$ in the above two cases, then the result follows.
\hfill\end{Proof}

\begin{figure}[htbp]
 \captionsetup{labelsep=period}
 \centering
 \begin{subfigure}[t]{0.35\linewidth}
  \centering
  \includegraphics[width=\linewidth]{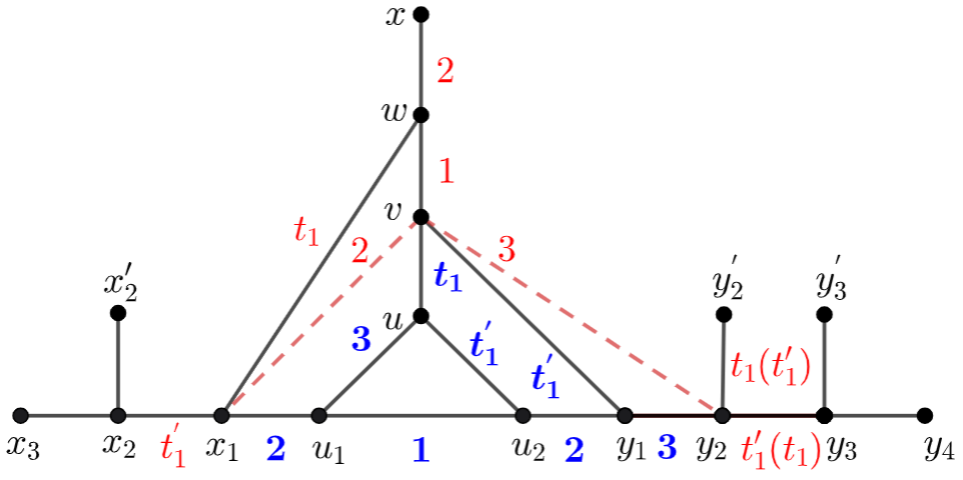}
  \caption{$3\notin\{\mu_2,\mu'_2\}$}
  \label{fig:ab1-a}
 \end{subfigure}
 \hspace{0.1\linewidth}
 \begin{subfigure}[t]{0.35\linewidth}
  \centering
  \includegraphics[width=\linewidth]{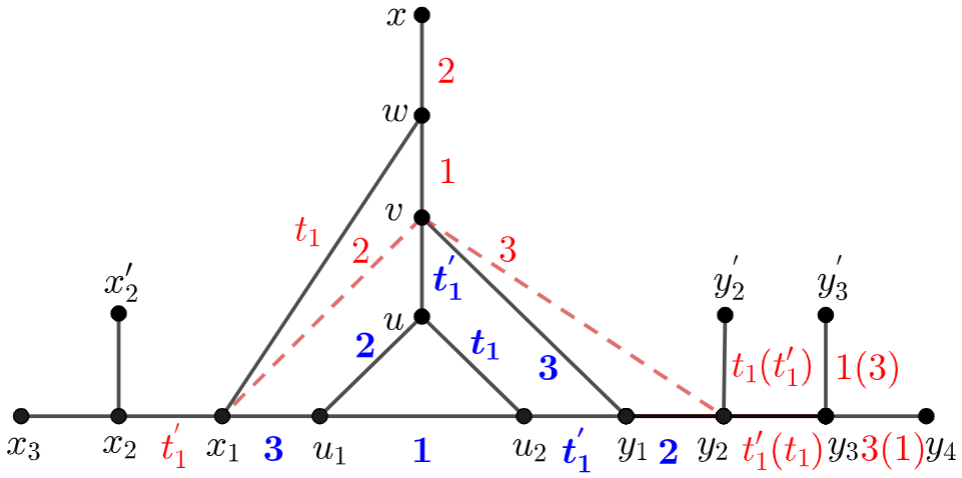}
  \caption{$3\notin\{\lambda_1,\lambda'_1\}$}
  \label{fig:ab1-b}
 \end{subfigure}
 \caption{a-b}
 \label{fig:ab1}
\end{figure}

\begin{figure}[htbp]
 \captionsetup{labelsep=period}
 \centering
 \begin{subfigure}[t]{0.31\linewidth}
  \centering
  \includegraphics[width=\linewidth]{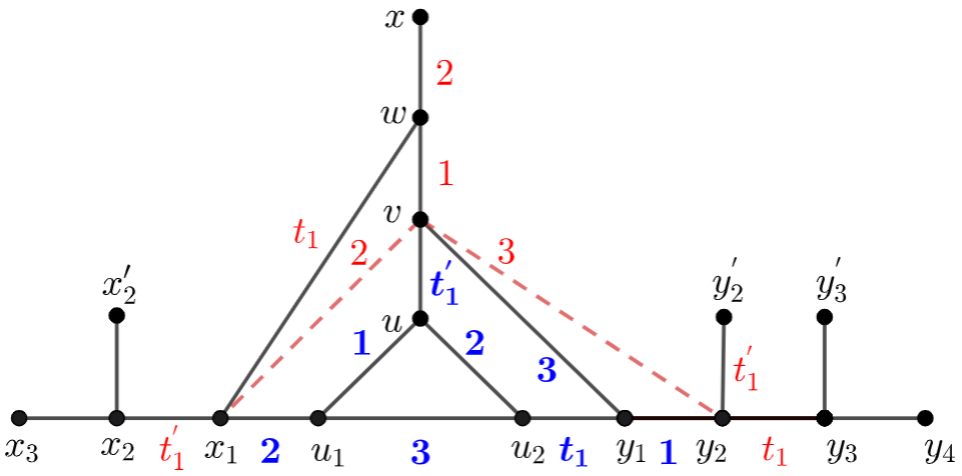}
  \caption{$\mu'_1=t_1$}
  \label{fig:ab2-a}
 \end{subfigure}
 \hspace{0.02\linewidth}
 \begin{subfigure}[t]{0.31\linewidth}
  \centering
  \includegraphics[width=\linewidth]{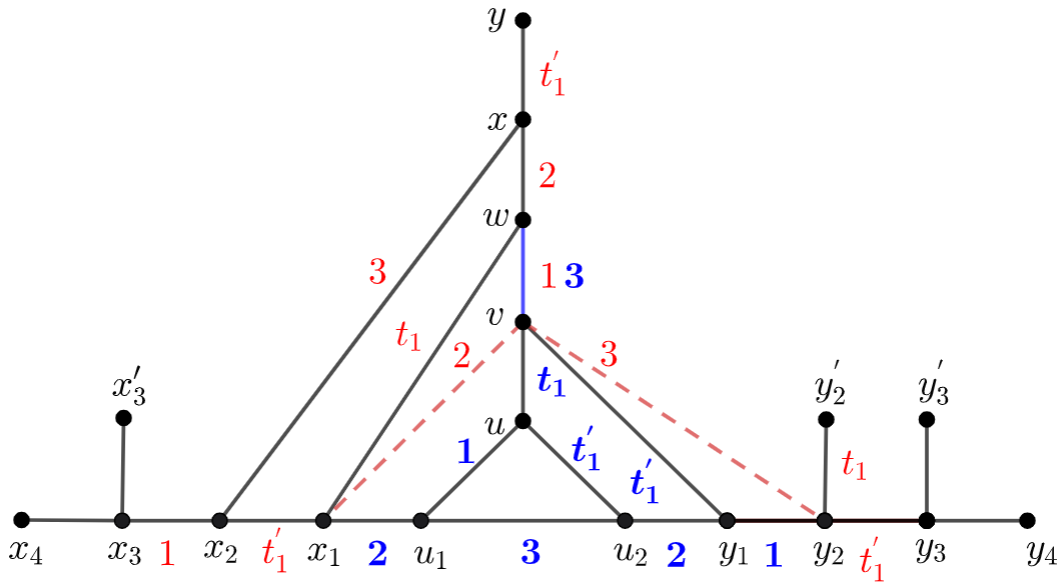}
  \caption{$\mu'_1=t'_1,\lambda'_1=1$}
  \label{fig:ab2-b}
 \end{subfigure}
 \hspace{0.02\linewidth}
 \begin{subfigure}[t]{0.31\linewidth}
  \centering
  \includegraphics[width=\linewidth]{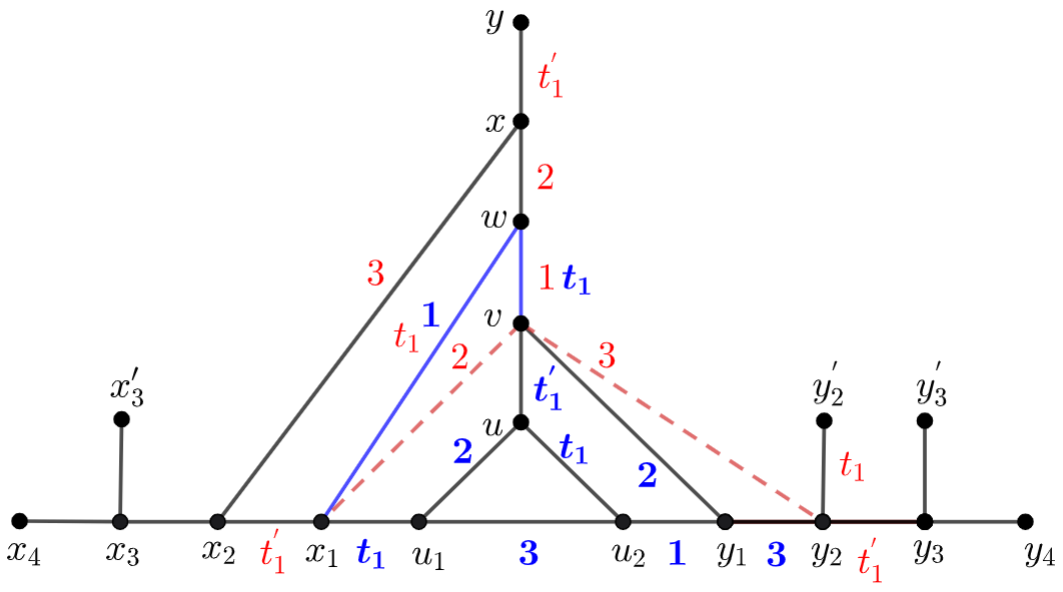}
  \caption{$\mu'_1=t'_1,\lambda'_1\neq1$}
  \label{fig:ab2-c}
 \end{subfigure}
 \caption{$1\notin\{\lambda_1,\lambda'_1\}$}
 \label{fig:abc2}
\end{figure}

\begin{figure}[htbp]
 \captionsetup{labelsep=period}
 \centering
 \begin{subfigure}[t]{0.35\linewidth}
  \centering
  \includegraphics[width=\linewidth]{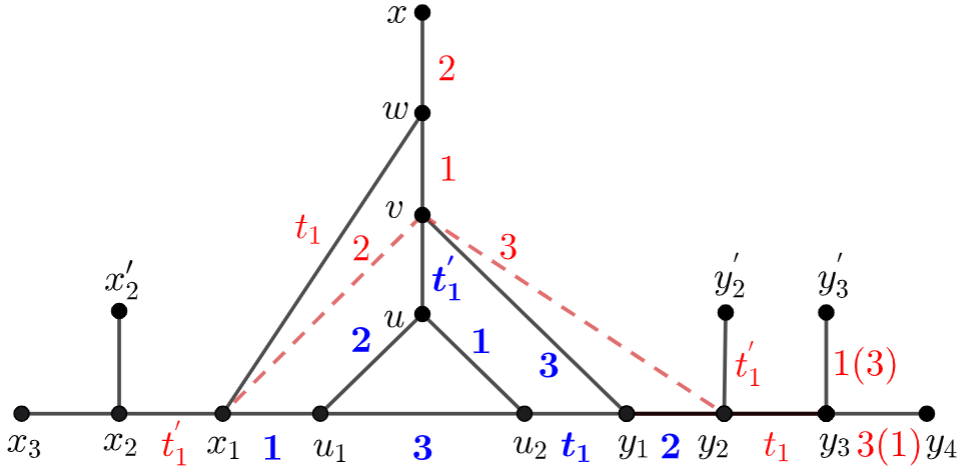}
  \caption{$\mu'_1=t_1$}
  \label{fig:ab3-a}
 \end{subfigure}
 \hspace{0.1\linewidth}
 \begin{subfigure}[t]{0.35\linewidth}
  \centering
  \includegraphics[width=\linewidth]{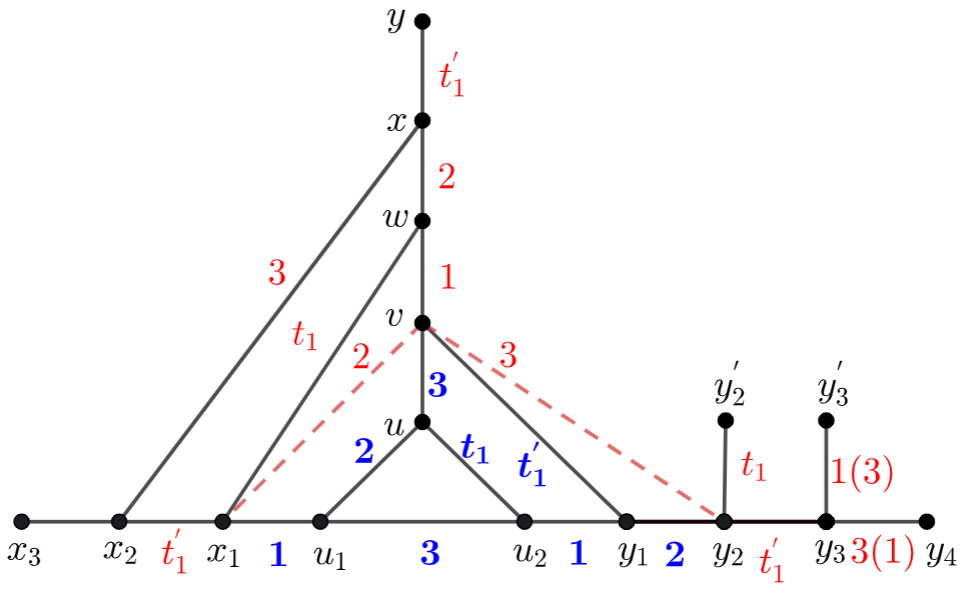}
  \caption{$\mu'_1=t'_1$}
  \label{fig:ab3-b}
 \end{subfigure}
 \caption{$1\notin\{\lambda_1,\lambda'_1\}$;~or~$\lambda'_1=1,t'_1\notin \{\lambda_2,\lambda'_2\}$}
 \label{fig:ab3}
\end{figure}

\begin{figure}[htbp]
 \captionsetup{labelsep=period}
 \centering
 \begin{subfigure}[t]{0.31\linewidth}
  \centering
  \includegraphics[width=\linewidth]{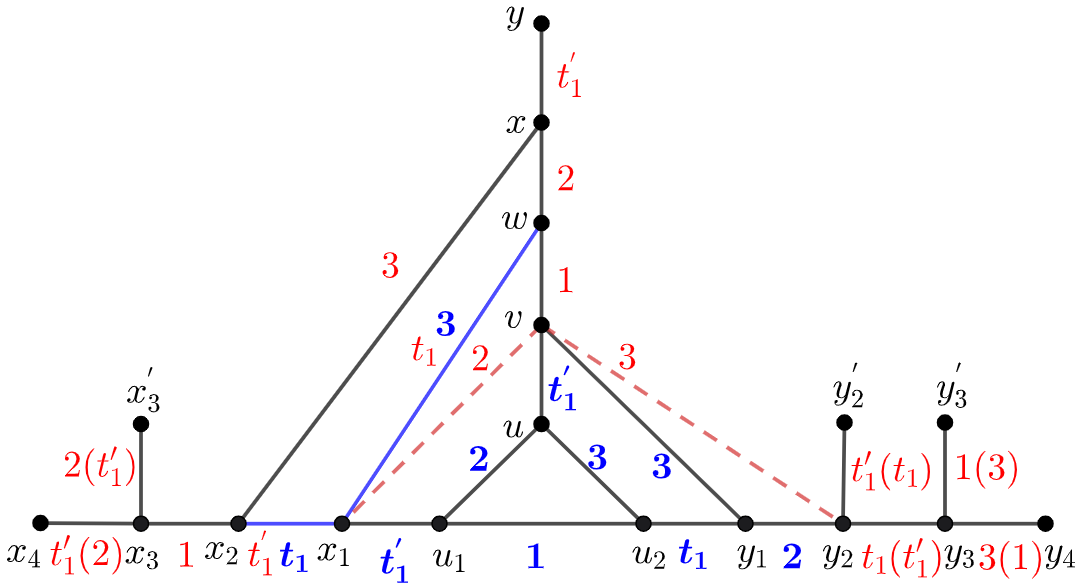}
  \caption{$\{\lambda_2,\lambda'_2\}=\{2,t'_1\}$}
  \label{fig:ab4-a}
 \end{subfigure}
 \hspace{0.02\linewidth}
 \begin{subfigure}[t]{0.31\linewidth}
  \centering
  \includegraphics[width=\linewidth]{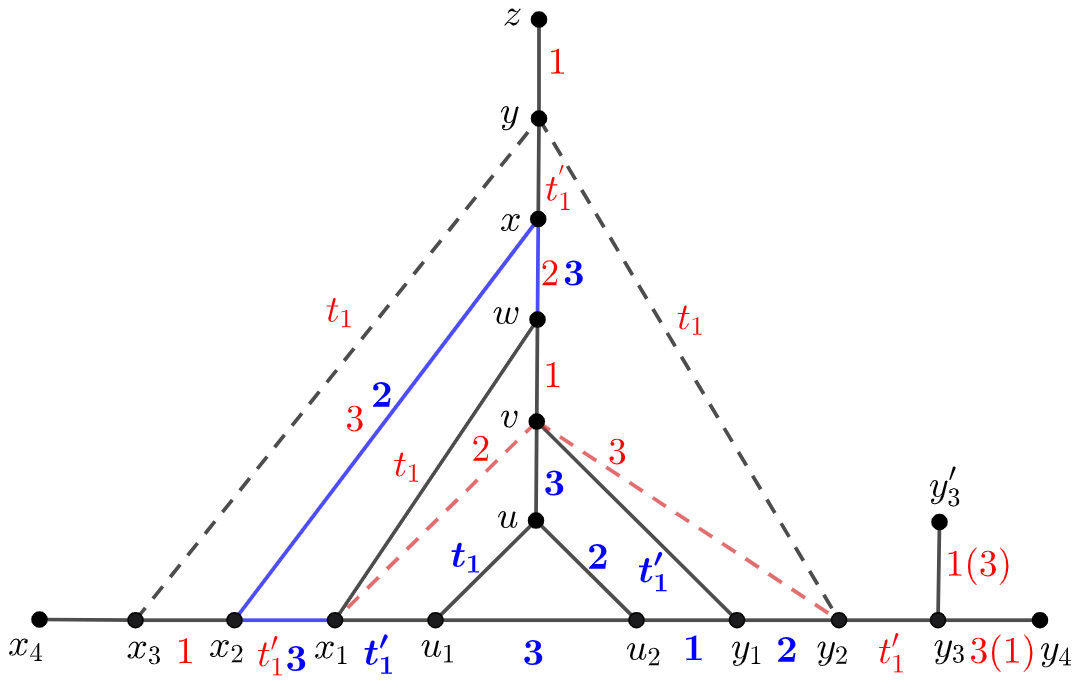}
  \caption{$2 \notin \{\lambda_2, \lambda'_2\},s_1=1$}
  \label{fig:ab4-b}
 \end{subfigure}
 \hspace{0.02\linewidth}
 \begin{subfigure}[t]{0.31\linewidth}
  \centering
  \includegraphics[width=\linewidth]{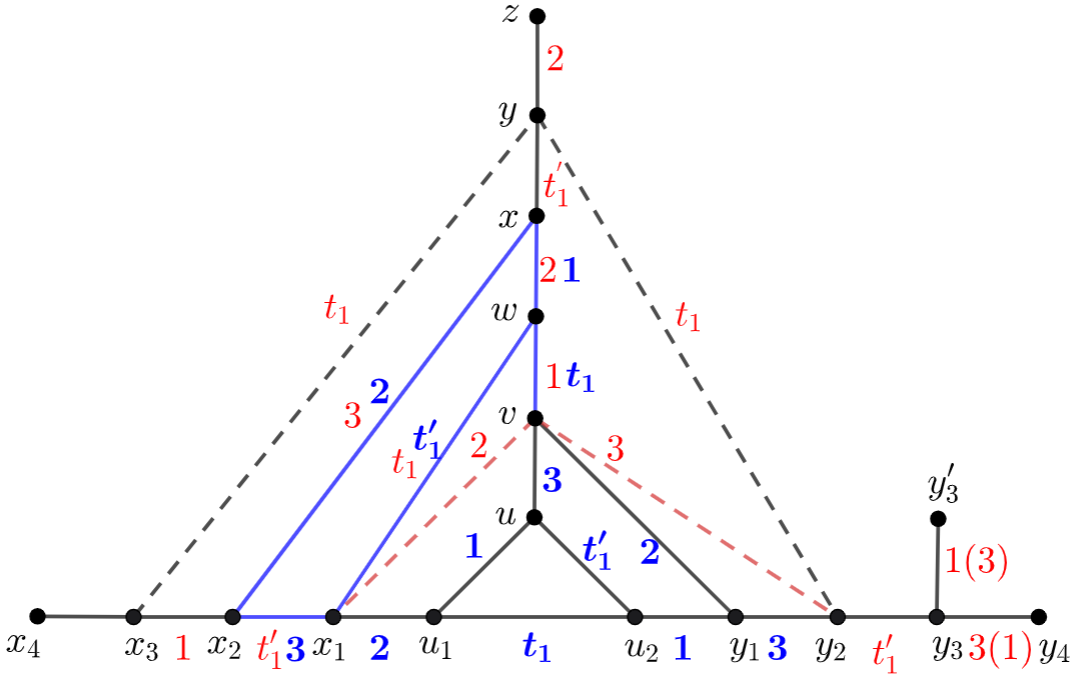}
  \caption{$2 \notin \{\lambda_2, \lambda'_2\},s_1=2$}
  \label{fig:ab4-c}
 \end{subfigure}
 \caption{$s_0=t'_1$}
 \label{fig:abc4}
\end{figure}

\begin{Corollary}\label{key2}
Let $N_{e_h}$ be a necklace with $h\neq2$, then $\chi'_{st}\left( N_{e_h} \right) =5$.
\end{Corollary}

\section{Conclusion}
\begin{figure}[htbp]
    \captionsetup{labelsep=period}
    \centering
    \begin{subfigure}[t]{0.3\linewidth}
        \centering
        \includegraphics[width=0.75\linewidth]{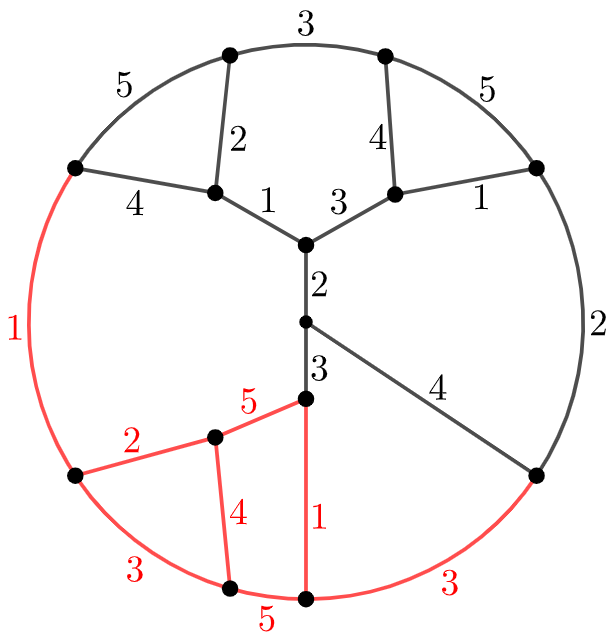}
        \caption{}
        \label{fig:ab-a}
    \end{subfigure}
    \hfill % Add horizontal space between subfigures
    \begin{subfigure}[t]{0.3\linewidth}
        \centering
        \includegraphics[width=0.75\linewidth]{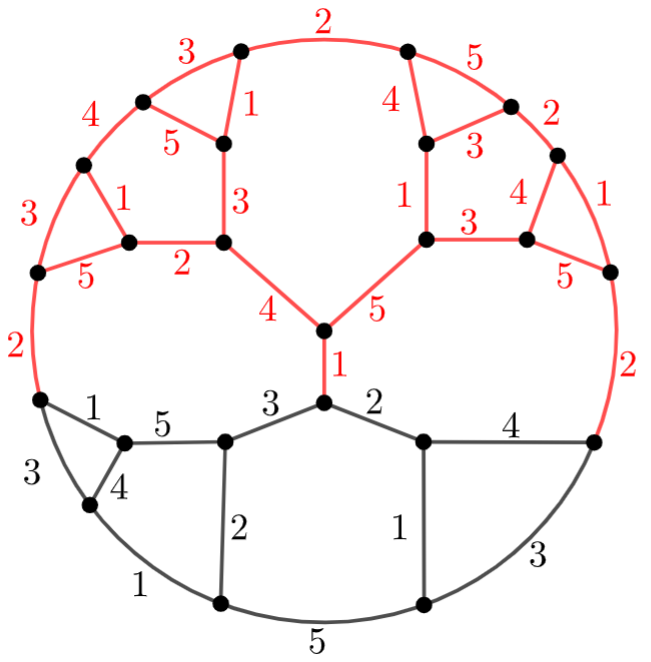}
        \caption{}
        \label{fig:ab-b}
    \end{subfigure}
    \hfill % Add horizontal space between subfigures
    \begin{subfigure}[t]{0.3\linewidth}
        \centering
        \includegraphics[width=0.75\linewidth]{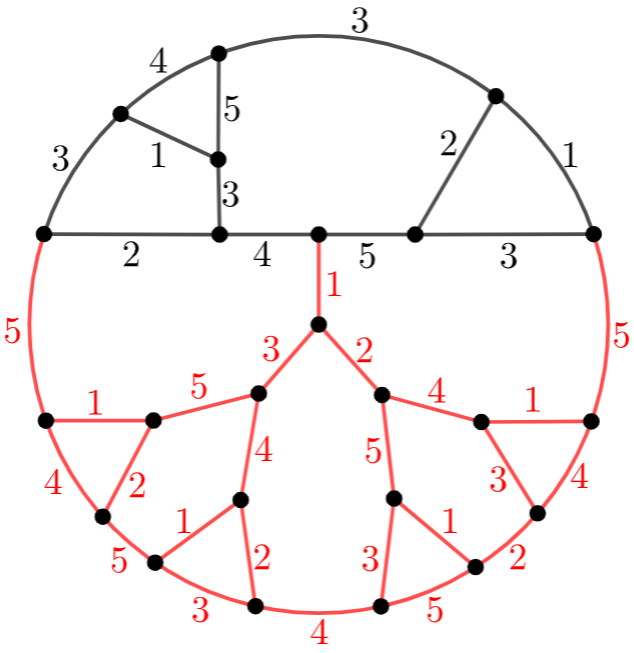}
        \caption{}
        \label{fig:ab-c}
    \end{subfigure}
    \caption{a-c}
    \label{fig10:abc}
\end{figure}

In this paper, we have determined star chromatic indices for the two  kinds of cubic Halin graphs. For the other kinds of cubic Halin graphs whose characteristic trees are neither complete graphs nor caterpillars, we believe that most of them obtain star chromatic index $5$, some of which can be deduced from Theorems~\ref{key3},~\ref{key4}, and~\ref{Th2}, such as the three graphs shown in Figure~\ref{fig10:abc}. So it is interesting to characterize which of such graphs have  star chromatic index $5$.

\section*{Author contributions}
Xingxing Hu: Investigation, Visualization, Writing-original draft; Yunfang Tang: Methodology,
Supervision, Writing-review and editing, Validation. All authors have read and approved the final
version of the manuscript for publication.

\section*{Use of Generative-AI tools declaration}
The authors declare that they have not used Artificial Intelligence (AI) tools in the creation of this
article.

%\section*{Acknowledgments}
%The authors are grateful to reviewers for their valuable comments and suggestions. This work
%is supported by the Natural Science Foundation of China (No.11701543).

\section*{Conflict of interest}
The authors declare there are no conflicts of interest.


\begin{thebibliography}{99}

\bibitem{BL2016} L. Bezegov\'a, B. Lu\v{z}ar, M. Mockov\v{c}iakov\'a, R. Sot\'ak, R. \v{S}krekovski, \textit{Star edge coloring of some classes of graphs}, J. Graph Theory., \textbf{81(1)} (2016), 73--82. https://doi.org/10.1002/jgt.21862

\bibitem{JCC2021} C. J. Casselgren, J. Granholm, A. Raspaud, \textit{On star edge colorings of bipartite and subcubic graphs}, Discrete Appl. Math., \textbf{298} (2021), 21--33. https://doi.org/10.1016/j.dam.2021.03.007

\bibitem{CL2012} G. J. Chang, D. Liu, \textit{Strong edge-coloring for cubic Halin graphs}, Discrete Math., \textbf{312(8)} (2012), 1468--1475. https://doi.org/10.1016/j.disc.2012.01.014

\bibitem{DMS2013} Z. Dvo\v{r}\'ak, B. Mohar, R. \v{S}\'amal, \textit{Star chromatic index}, J. Graph Theory., \textbf{72(3)} (2013), 313--326. https://doi.org/10.1002/jgt.21644

\bibitem{Hou2019} X. Hou, L. Li, T. Wang, \textit{Star edge coloring of some special graphs}, arXiv:2010.14349 [math.CO], (2020). Unpublished work. https://doi.org/10.48550/arXiv.2010.14349

\bibitem{LD2008} X. Liu, K. Deng, \textit{An upper bound on the star chromatic index of graphs with $\Delta\ge 7$}, J. Lanzhou Univ. (Nat. Sci.), \textbf{44} (2008), 98--99. https://doi.org/10.13885/j.issn.0455-2059.2008.02.003

\bibitem{LL2012} K. W. Lih, D. Liu, \textit{On the strong chromatic index of cubic Halin graphs}, Appl. Math. Lett., \textbf{25} (2012), 898--901. https://doi.org/10.1016/j.aml.2011.10.046

\bibitem{lei2018} H. Lei, Y. Shi, Z. X. Song, \textit{Star chromatic index of subcubic multigraphs}, J. Graph Theory., \textbf{88(4)} (2018), 566--576. https://doi.org/10.1002/jgt.22230

\bibitem{lei2018b} H. Lei, Y. Shi, Z. X. Song, T. Wang, \textit{Star 5-edge-colorings of subcubic multigraphs}, Discrete Math., \textbf{341(4)} (2018), 950--956. https://doi.org/10.1016/j.disc.2017.12.008

\bibitem{LS2021} H. Lei, Y. Shi, \textit{A survey on star edge colorings of graphs}, Prog. Math., (China) \textbf{50(1)} (2021), 77--93. https://doi.org/10.11845/sxjz.2020006a

\bibitem{Shiu2006} W. C. Shiu, P. C. B. Lam, W. K. Tam, \textit{On strong chromatic index of Halin graphs}, Congr. Numer., \textbf{57} (2006), 211--222. https://www.researchgate.net/publication/238636258

\bibitem{ShiuT2009} W. C. Shiu, W. K. Tam, \textit{The strong chromatic index of complete cubic Halin graphs}, Appl. Math. Lett., \textbf{22} (2009), 754--758. https://doi.org/10.1016/j.aml.2008.08.019

\bibitem{WY2021} W. Yang, B. Wu, \textit{Proof of a conjecture on the strong chromatic index of Halin graphs}, Discrete Appl. Math., \textbf{302} (2021), 92--102. https://doi.org/10.1016/j.dam.2021.06.016

\end{thebibliography}
\end{document}